# The quenched invariance principle for random walks in random environments admitting a bounded cycle representation

Jean-Dominique Deuschel[a]  and Holger Kösters[b]

[a]*Fachbereich Mathematik, Sekr. MA 7–4, Technische Universität Berlin, Straße des 17. Juni 136, D-10623 Berlin, Germany.*
*E-mail: deuschel@math.tu-berlin.de*
[b]*Fakultät für Mathematik, Universität Bielefeld, Postfach 100131, D-33501 Bielefeld, Germany.*
*E-mail: hkoesters@math.uni-bielefeld.de*

**Abstract.** We derive a quenched invariance principle for random walks in random environments whose transition probabilities are defined in terms of weighted cycles of bounded length. To this end, we adapt the proof for random walks among random conductances by Sidoravicius and Sznitman (*Probab. Theory Related Fields* **129** (2004) 219–244) to the non-reversible setting.

**Résumé.** Nous dérivons un principe d'invariance presque sûr pour les marches aléatoires en milieu aléatoire dont les transitions sont données par des poids indexés par des cycles bornés. A cet effet nous adaptons la démonstration pour les marches symétriques en milieu aléatoire de Sidoravicius et Sznitman (*Probab. Theory Related Fields* **129** (2004) 219–244) dans le cas non réversible.

## 1. Introduction and summary

We consider a class of random walks in random environments (RWRE's) on $\mathbb{Z}^d$ admitting the following "bounded cycle representation":

We begin by introducing some terminology. A *cycle* $C$ is a finite sequence $(z_0, z_1, \ldots, z_n)$ of points in $\mathbb{Z}^d$ such that $z_0, \ldots, z_{n-1}$ are pairwise different and $z_n = z_0$. The number $n$ is also called the *length* of the cycle. We allow cycles of lengths 1 and 2. For $x \in \mathbb{Z}^d$, we write $x \in C$ if there is an index $i = 0, \ldots, n-1$ such that $z_i = x$. For $x \in \mathbb{Z}^d$, the cycle $C + x$ is defined by $(z_0 + x, z_1 + x, \ldots, z_n + x)$. Finally, we also identify the cycle $C$ with the sequence of its (directed) edges $((z_0, z_1), (z_1, z_2), \ldots, (z_{n-1}, z_n))$. Thus, for $(x, y) \in \mathbb{Z}^d \times \mathbb{Z}^d$, we write $(x, y) \in C$ if there is an index $i = 0, \ldots, n-1$ such that $z_i = x$, $z_{i+1} = y$.

Let $K \in \mathbb{N}$, and let $C_1, \ldots, C_K$ be cycles of lengths $n_1, \ldots, n_K$. It may be helpful to think of all these cycles having a common point (e.g. the origin), although we do not need that. Let $(\Omega, \mathcal{F}, \mathbb{P})$ be a probability space endowed with a group of measurable transformations $(T_x)_{x \in \mathbb{Z}^d}$ such that $T_{x+y} = T_x \circ T_y$ for all $x, y \in \mathbb{Z}^d$, let $W_1, \ldots, W_K$ be non-negative random variables on $(\Omega, \mathcal{F}, \mathbb{P})$, and let $M$ be the random variable on $(\Omega, \mathcal{F}, \mathbb{P})$ given by

$$M(\omega) := \sum_{i=1}^{K} \sum_{x \in \mathbb{Z}^d} W_i(T_x \omega) \cdot \mathbf{1}_{\{0 \in C_i + x\}} \tag{1.1}$$

for all $\omega \in \Omega$. We suppose that $M$ is strictly positive. (See also Assumption (b) below.) For each $\omega \in \Omega$ and each $z \in \mathbb{Z}^d$, set

$$p_z(\omega) := \frac{1}{M(\omega)} \cdot \sum_{i=1}^{K} \sum_{x \in \mathbb{Z}^d} W_i(T_x \omega) \cdot \mathbf{1}_{\{(0,z) \in C_i + x\}}. \tag{1.2}$$

Then, for any $\omega \in \Omega$, we have $p_z(\omega) \geq 0$ for all $z \in \mathbb{Z}^d$ and $\sum_{z \in \mathbb{Z}^d} p_z(\omega) = 1$, i.e. the $p_z(\omega)$ define a probability measure on $\mathbb{Z}^d$. For fixed $x \in \mathbb{Z}^d$ and fixed $\omega \in \Omega$, the random walk in the "environment" $\omega$ and with the "starting point" $x$ is the Markov chain $(X_n)_{n \in \mathbb{N}}$ on the probability space $(\Sigma, \mathcal{G}, P_{x,\omega})$ with transition probabilities

$$P_{x,\omega}(X_{n+1} = y + z | X_n = y) = p_z(T_y \omega)$$

and initial distribution

$$P_{x,\omega}(X_0 = x) = 1.$$

Obviously, there is no loss of generality in assuming that $(\Sigma, \mathcal{G})$ is the set $(\mathbb{Z}^d)^\mathbb{N}$ endowed with the product $\sigma$-field, $(X_n)_{n \in \mathbb{N}}$ is the coordinate process, and $P_{x,\omega}$ is the probability measure on $(\Sigma, \mathcal{G})$ determined by the above conditions.

First of all, note that the RWRE has *bounded range*, uniformly in $\omega$, because

$$\Lambda := \{z \in \mathbb{Z}^d \mid \exists i = 1, \ldots, K, x \in \mathbb{Z}^d: (0, z) \in C_i + x\} \tag{1.3}$$

is a finite set and for any $\omega \in \Omega$, we have

$$p_z(\omega) = 0 \quad \text{for all } z \notin \Lambda. \tag{1.4}$$

We continue by stating our standing assumptions. In the sequel let $e_i$ denote the $i$th unit vector in $\mathbb{Z}^d$ ($i = 1, \ldots, d$), and let $\|\cdot\|_2$ denote the Euclidean norm on $\mathbb{R}^d$.

*Assumptions.*

(a) *The probability measure $\mathbb{P}$ is invariant and ergodic with respect to each of the subgroups $(T_{ze_i})_{z \in \mathbb{Z}}$ ($i = 1, \ldots, d$).*
(b) *The random variable $M$ is bounded away from $0$ and $\infty$, i.e. there exist constants $0 < c \leq C < \infty$ such that $c \leq M(\omega) \leq C$ for all $\omega \in \Omega$.*
(c) *The RWRE is strongly irreducible, uniformly in $\omega$, i.e. there exist $\varepsilon > 0$ and $N \in \mathbb{N}$ such that for all $x \in \mathbb{Z}^d$, for all $\omega \in \Omega$ and for all $e \in \mathbb{Z}^d$ with $\|e\|_2 = 1$, we have $P_{x,\omega}(X_n = x + e) \geq \varepsilon$ for some $n \leq N$.*

*Remarks.*

Observe that Assumption (a) entails the following weaker condition:
(a′) *The probability measure $\mathbb{P}$ is invariant and ergodic with respect to the group $(T_x)_{x \in \mathbb{Z}^d}$.*

Also, due to (1.4), Assumption (c) is equivalent to the following condition:
(c′) *There exist $\varepsilon_0 > 0$ and $N \in \mathbb{N}$ such that for all $x \in \mathbb{Z}^d$, for all $\omega \in \Omega$ and for all $e \in \mathbb{Z}^d$ with $\|e\|_2 = 1$, there exist $n \leq N$ and $z_0, \ldots, z_n \in \mathbb{Z}^d$ with $z_0 = x$, $z_n = x + e$ and $p_{z_i - z_{i-1}}(T_{z_{i-1}} \omega) \geq \varepsilon_0$ for all $i = 1, \ldots, n$.*

In particular, this holds under the following stronger condition that the RWRE is *uniformly elliptic*:

(c″) *The set $\Lambda$ defined in (1.3) generates $\mathbb{Z}^d$, and there exists a constant $\varepsilon_0 > 0$ such that $\inf_{z \in \Lambda} p_z(\omega) \geq \varepsilon_0$ for all $\omega \in \Omega$.*

However, see Example (c) below.

***Examples.***

(a) The random conductance model. *Take $K = d$ and $C_i = (0, e_i, 0)$ $(i = 1, \ldots, K)$, and suppose that the random variables $W_i$ take values in $[a, b]$ $(i = 1, \ldots, K)$ for some $0 < a < b < \infty$. Then Assumptions (b) and (c″) are clearly satisfied.*

(b) The uniformly elliptic case. *Suppose that the cycles $C_i$ are such that $\Lambda$ generates $\mathbb{Z}^d$ and that the random variables $W_i$ take values in $[a, b]$ $(i = 1, \ldots, K)$ for some $0 < a < b < \infty$. Then Assumptions (b) and (c″) are clearly satisfied. For the case where the random variable $M$ is constant, this model has been introduced in Section 4.3 in [7] (see also pp. 124–125 in [8]).*

(c) The square triangle model. *Take $d = 2$, $K = 2$, $C_1 = ((0,0), (1,0), (1,1), (0,1), (0,0))$, $C_2 = ((0,0), (1,0), (1,1), (0,0))$, and suppose that the random variables $W_i$ take values in $\{0, 1\}$ such that $W_1(\omega) + W_2(\omega) = 1$ for all $\omega \in \Omega$ and that the random vectors $(W_1 \circ T_x, W_2 \circ T_x)$, $x \in \mathbb{Z}^2$, are i.i.d. It is obvious that this model satisfies Assumptions (b) and (c) (with $N = 2$), but not Assumption (c″).*

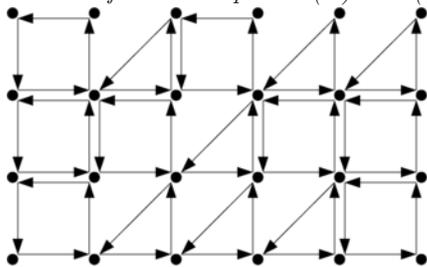

(d) The triangle triangle model. *Take $d = 2$, $K = 2$, $C_1 = ((0,0), (1,1), (0,1), (0,0))$, $C_2 = ((0,0), (1,0), (1,1), (0,0))$, and suppose that the random variables $W_i$ take values in $\{0, 1\}$ such that $W_1(\omega) + W_2(\omega) = 1$ for all $\omega \in \Omega$ and that the random vectors $(W_1 \circ T_x, W_2 \circ T_x)$, $x \in \mathbb{Z}^2$, are i.i.d. It is obvious that this model satisfies Assumption (b), but not Assumption (c), as it contains "corridors" of arbitrary length. Thus it does not fit into our framework.*

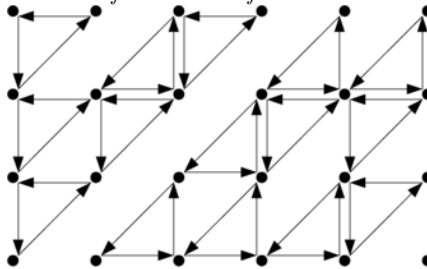

(e) The random walk on the supercritical percolation cluster. *Take $K = d$ and $C_i = (0, e_i, 0)$ $(i = 1, \ldots, K)$, suppose that the random variables $W_i$ are independent with $P(W_i = 1) = p = 1 - P(W_i = 0)$ for some $p \in (p_{\text{crit}}, 1)$, where $p_{\text{crit}}$ denotes the critical percolation probability, and condition upon the event that the origin belongs to the unique infinite open cluster. Clearly, this model does not fit into our framework either.*

In the study of the random walk $(X_n)_{n \in \mathbb{N}}$, one usually distinguishes between "quenched" results relating to the probability measures $P_{0,\omega}$, for any $\omega \in \Omega$, and "annealed" results relating to the averaged probability measure given by

$$(\mathbb{P} \circ P_{0,\cdot})(A) := \int_\Omega P_{0,\omega}(A) \, \mathbb{P}(d\omega), \quad A \in \mathcal{G}.$$

Quenched results are usually harder to prove, but they are also more relevant.

We are interested in the question whether the RWRE satisfies the quenched Invariance Principle. More precisely, for $N \geq 1$, set

$$\beta_N := \text{the polygonal interpolation of } \frac{k}{N} \mapsto \frac{X_k}{\sqrt{N}}, \quad k \geq 0,$$

and note that $\beta_N$ is a random variable taking values in $\mathscr{C}([0,\infty[;\mathbb{R}^d)$, the space of continuous functions on $[0,\infty[$ taking values in $\mathbb{R}^d$ endowed with the topology of uniform convergence on compact intervals. The quenched Invariance Principle states that for $\mathbb{P}$-a.e. $\omega \in \Omega$, under the probability measure $P_{0,\omega}$, the random variables $\beta_N$ converge in distribution to a $d$-dimensional Brownian motion with a zero mean and a non-trivial covariance matrix.

For model (a), the quenched Invariance Principle has recently been established by Sidoravicius and Sznitman (see Theorem 1.1 and Remark 1.3 in [14]). The main aim of this paper is to extend their result to the more general case of a RWRE admitting a bounded cycle representation:

**Theorem 1.1.** *Suppose that the RWRE admits a bounded cycle representation. Then, for $\mathbb{P}$-a.e. $\omega \in \Omega$, under the probability measure $P_{0,\omega}$, the random variables $\beta_N$ converge in distribution to a $d$-dimensional Brownian motion with mean zero and non-degenerate covariance matrix $A$ as in* (4.1) *below.*

For model (b), in the special case where the random variable $M$ is constant, it follows from a recent result by Komorowski and Olla [7] that (already under the weaker assumption (a$'$)) the RWRE satisfies the annealed Central Limit Theorem, i.e. under the probability measure $\mathbb{P} \circ P_{0,.}$, the random vectors $Z_N := X_N/\sqrt{N}$ converge in distribution to a $d$-dimensional Gaussian random variable with a zero mean and a non-trivial covariance matrix (see Theorem 2.2 and Example 4.3 in [7]). The above theorem shows that (under the stronger assumption (a)) the RWRE also satisfies the quenched Central Limit Theorem, i.e. for $\mathbb{P}$-a.e. $\omega \in \Omega$, under the probability measure $P_{0,\omega}$, the random vectors $Z_N := X_N/\sqrt{N}$ converge in distribution to a $d$-dimensional Gaussian random variable with a zero mean and a non-trivial covariance matrix.

Furthermore, for model (e), the quenched invariance principle has recently been proved by Sidoravicius and Sznitman [14] in dimension $d \geq 4$, and by Berger and Biskup [1], Mathieu and Piatnitski [13] in dimension $d \geq 2$. Here one can take advantage of very precise results on the transition kernel of the random walk on the infinite percolation cluster.

Finally, let it be mentioned that Mathieu [12] and Biskup and Prescott [2] have very recently proved the quenched invariance principle for certain variants of model (a) in which the random variables $W_i$ need not be bounded away from 0.

In proving the above theorem, we will closely follow the approach of Sidoravicius and Sznitman [14] for the random conductance model (see Section 1 in [14]). To this end, we need two basic ingredients:

- We introduce the so-called *corrector function* $\chi: \mathbb{Z}^d \times \Omega \longrightarrow \mathbb{R}^d$, for which the process $(M_n^\omega)_{n \in \mathbb{N}}$ defined by

  $$M_n^\omega := X_n + \chi(X_n, \omega), \quad n \in \mathbb{N},$$

  is a martingale under $P_{0,\omega}$, for $\mathbb{P}$-almost each $\omega \in \Omega$.

- We establish the upper Gaussian bound

  $$P_{x,\omega}(X_n = y) \leq C_0 n^{-d/2} \exp\left(-\frac{\|x-y\|_2^2}{C_0 n}\right), \quad n \geq 1, \; x,y \in \mathbb{Z}^d,$$

  on the transition kernel of the random walk under $P_{0,\omega}$, for each $\omega \in \Omega$. Here $C_0 > 0$ is a constant not depending on $\omega$.

In Sidoravicius and Sznitman [14] (and also in Berger and Biskup [1] and Mathieu and Piatnitski [13]), the construction of the corrector function relies on the reversibility of the so-called chain of environments viewed from the particle (see Section 3 for further details). Unfortunately, this property is generally lost in case that the RWRE admits a bounded cycle representation when one of the representing cycles $C_i$ has length $n_i \geq 3$. We therefore have to take a somewhat different approach when constructing the corrector function, using ideas from the theory of (asymmetric) exclusion processes (see e.g. [9, 10, 15]).

We have the following upper Gaussian bound on the transition kernel of the random walk:

**Proposition 1.2.** *Suppose that the RWRE admits a bounded cycle representation. There exists a constant $C_0 > 0$ (not depending on $\omega$) such that for all $n \geq 1$ and all $x, y \in \mathbb{Z}^d$,*

$$P_{x,\omega}(X_n = y) \leq C_0 n^{-d/2} \exp\left(-\frac{\|y-x\|_2^2}{C_0 n}\right).$$

Such upper Gaussian bounds are well known for reversible random walks (see e.g. Section 14 in [16]), but they generally fail for non-reversible random walks. Now a RWRE admitting a bounded cycle representation is generally not reversible when one of the representing cycles $C_i$ has length $n_i \geq 3$, but it belongs to the class of so-called centered random walks (see Section 2 for further details). For this class, it has recently been shown by Mathieu [11] that the "off-diagonal" upper bound on the transition kernel follows from the "on-diagonal" upper bound on the transition kernel, which will be obtained by means of a Nash inequality. However, also in this context, there are some complications arising from the lack of reversibility.

This paper is structured as follows. In Section 2 we give the proof of Proposition 1.2. In Section 3 we construct the corrector function. After these preparations, the proof of Theorem 1.1 is virtually identical to that of Theorem 1.1 in [14]; Section 4 contains some relevant comments. Finally, Appendices A and B contain a number of auxiliary results which seem either quite standard or very similar to existing results and which have been included for the sake of completeness.

## 2. The upper Gaussian bound on the transition kernel of the random walk

A *centered random walk* as introduced in Definition 2.1 in [11] is an irreducible Markov chain $(X_n)_{n \in \mathbb{N}}$ taking values in $\mathbb{Z}^d$ for which there exist a collection of cycles $(C_i)_{i \in \mathbb{N}}$ of bounded length, a collection of weights $(w_i)_{i \in \mathbb{N}}$ and a positive measure $\pi$ on $\mathbb{Z}^d$ (a so-called *centering measure*) such that for all $n \in \mathbb{N}$ and for all $x, y \in \mathbb{Z}^d$,

$$P(X_{n+1} = y | X_n = x) = \frac{1}{\pi(\{x\})} \cdot \sum_{i \in \mathbb{N}} w_i \cdot \mathbf{1}_{\{(x,y) \in C_i\}}.$$

It is immediate from our definitions that for each $\omega \in \Omega$, the random walk $(X_n)_{n \in \mathbb{N}}$ under $P_{0,\omega}$ is a centered random walk in this sense. More precisely, the measure on $\mathbb{Z}^d$ given by

$$\pi_\omega(\{x\}) := M(T_x \omega), \quad x \in \mathbb{Z}^d,$$

is a centering measure, and hence also an invariant measure (see Lemma 2.5 in [11]), for the random walk $(X_n)_{n \in \mathbb{N}}$ under $P_{0,\omega}$. Moreover, by our Assumption (b), we have

$$\inf_{x \in \mathbb{Z}^d} \pi_\omega(\{x\}) \geq \inf_{\omega \in \Omega} M(\omega) \geq c > 0, \tag{2.1}$$

as required for most of the results in [11].

As already mentioned, the Markov chain $(X_n)_{n \in \mathbb{N}}$ under $P_{0,\omega}$ is generally *not* reversible (with respect to $\pi_\omega$) when one of the representing cycles $C_i$ has length $n_i \geq 3$. At least, we have an explicit description of the time-reversed Markov chain in this case: It is given by the "reversed" cycles.

More precisely, for any cycle $C = (x_0, x_1, \ldots, x_n)$, the reversed cycle is given by $C^* = (x_n, x_{n-1}, \ldots, x_0)$. Let $C_1^*, \ldots, C_K^*$ be the reversed cycles belonging to the cycles $C_1, \ldots, C_K$, and let $W_1, \ldots, W_K$ and $M$ be the same as before. Then we clearly have

$$M(\omega) = \sum_{i=1}^{K} \sum_{x \in \mathbb{Z}^d} W_i(T_x \omega) \cdot \mathbf{1}_{\{0 \in C_i^* + x\}}$$

for any $\omega \in \Omega$, so that we may introduce the probabilities $p_z^*(\omega)$ (defined analogously to the probabilities $p_z(\omega)$) and the probability measures $P_{0,\omega}^*$ (defined analogously to the probability measures $P_{0,\omega}$). As observed in Section 2.2 in [11], the Markov chain $(X_n)_{n\in\mathbb{N}}$ under $P_{0,\omega}^*$ is then linked to the time-reversed Markov chain of the Markov chain $(X_n)_{n\in\mathbb{N}}$ under $P_{0,\omega}$. That is, it also has $\pi_\omega$ as an invariant measure, and

$$p_z^*(T_y\omega) = \frac{M(T_{y+z}\omega)}{M(T_y\omega)} \cdot p_{-z}(T_{y+z}\omega) \tag{2.2}$$

for all $y, z \in \mathbb{Z}^d$.

In particular, the time-reversed Markov chain is of the same type as the original Markov chain. This observation plays a crucial role in the derivation of the upper Gaussian bound in Theorem 2.10 in [11], and also in the derivation of the bound (2.3).

We now explain how to prove Proposition 1.2. Suppose that there is a constant $C_1 > 0$ (not depending on $\omega$) such that

$$P_{x,\omega}(X_n = y) \leq C_1 M(T_y\omega) n^{-d/2} \tag{2.3}$$

for all $n \geq 1$ and all $x, y \in \mathbb{Z}^d$. Then, by Theorem 2.10 in [11], there is a constant $C_2 > 0$ depending only on $c$ (from (2.1)), $\max\{n_1, \ldots, n_K\}$ (from Section 1) and $C_1$ such that

$$P_{x,\omega}(X_n = y) \leq C_2 M(T_y\omega) n^{-d/2} \exp\left(-\frac{d(x,y)^2}{C_2 n}\right) \tag{2.4}$$

for all $n \geq 1$ and all $x, y \in \mathbb{Z}^d$. Here, $d(x,y)$ denotes the natural graph distance associated with the random walk $(X_n)_{n\in\mathbb{N}}$ under the probability measure $P_{0,\omega}$ (see Section 2.1 in [11]). Since the latter has bounded range $B$ (with respect to the Euclidean distance), uniformly in $\omega \in \Omega$, it is related to the Euclidean distance by the inequality $d(x,y) \geq \|y - x\|_2 / B$. Also, by our Assumption (b), the factor $M(T_y\omega)$ is bounded above by the constant $C$, uniformly in $\omega \in \Omega$. Thus, there is a constant $C_3 > 0$ depending only on $B$, $C$ and $C_2$ such that

$$P_{x,\omega}(X_n = y) \leq C_3 n^{-d/2} \exp\left(-\frac{\|y - x\|_2^2}{C_3 n}\right) \tag{2.5}$$

for all $n \geq 1$ and all $x, y \in \mathbb{Z}^d$. In particular, $C_3$ is independent of $\omega$.

To complete the proof of Proposition 1.2, it remains to verify (2.3). This will be done by means of a Nash inequality. To this end, we introduce some more notation. Given a transition kernel $Q$ on $\mathbb{Z}^d$ with (positive) invariant measure $\pi$, let

$$\mathcal{E}_Q(f,f) := \frac{1}{2}\sum_{x,y}(f(y) - f(x))^2 \pi(x) Q(x,y) = \langle f, (I - Q)f\rangle_\pi$$

denote the associated *Dirichlet form*, where $\langle \cdot, \cdot \rangle_\pi$ denotes the scalar product in $L^2(\pi)$. Furthermore, let $Q^*$ denote the adjoint of $Q$ with respect to $\pi$, i.e. let

$$Q^*(x,y) := \frac{\pi(\{y\})}{\pi(\{x\})} \cdot Q(y,x)$$

for any $x, y \in \mathbb{Z}^d$. Finally, let $\ell_0(\mathbb{Z}^d)$ denote the set of functions on $\mathbb{Z}^d$ with finite support.

**Proposition 2.1.** *Let $(X_n)_{n\in\mathbb{N}}$ be a random walk on $\mathbb{Z}^d$ with transition kernel $Q$ and invariant measure $\pi$ for which there exist constants $0 < c \leq C < \infty$ with $c \leq \pi(\{x\}) \leq C$ for all $x \in \mathbb{Z}^d$. Suppose that the transition kernel $Q$ is strongly irreducible, i.e.*

$$\exists \varepsilon > 0, \ \exists N \in \mathbb{N}, \ \forall x \in \mathbb{Z}^d, \ \forall e \in \mathbb{Z}^d \colon \|e\|_2 = 1, \ \exists n \leq N, \quad Q^n(x, x+e) \geq \varepsilon,$$

*and of bounded range, i.e.*

$$\exists B > 0, \ \forall x, y \in \mathbb{Z}^d, \quad \|y - x\|_2 > B \implies Q(x, y) = 0.$$

*Then there exists a number $m \geq 1$ such that the transition kernel $(Q^m)^* Q^m$ satisfies the d-dimensional Nash inequality, i.e.*

$$\exists \kappa > 0, \ \forall f \in \ell_0(\mathbb{Z}^d), \quad \mathcal{E}_{(Q^m)^* Q^m}(f, f) \geq \kappa \|f\|_{L^2(\pi)}^{2+4/d} \|f\|_{L^1(\pi)}^{-4/d}.$$

*More precisely, $m$ and $\kappa$ are constants depending only on $c$, $C$, $\varepsilon$, $N$ and $B$.*

Although we think that this result should be well known, we are not aware of a suitable reference in the literature. (Most of the existing accounts of Nash inequalities seem to concentrate on reversible Markov chains.) A full proof of Proposition 2.1 is therefore given in Appendix A of this paper.

Besides that, we will need the following result (see also Theorem 4.1 in [4]):

**Proposition 2.2.** *Let $(X_n)_{n \in \mathbb{N}}$ be a random walk on $\mathbb{Z}^d$ with transition kernel $Q$ and invariant measure $\pi$, for which there exist constants $0 < c \leq C < \infty$ with $c \leq \pi(\{x\}) \leq C$ for all $x \in \mathbb{Z}^d$. Suppose that the transition kernel $Q^* Q$ satisfies the d-dimensional Nash inequality, i.e.*

$$\exists \kappa > 0, \ \forall f \in \ell_0(\mathbb{Z}^d), \quad \mathcal{E}_{Q^* Q}(f, f) \geq \kappa \|f\|_{L^2(\pi)}^{2+4/d} \|f\|_{L^1(\pi)}^{-4/d}.$$

*Then there exists a constant $C_0$ (depending only on $c$, $C$ and $\kappa$) such that*

$$\|Q^n\|_{L^1(\pi) \to L^2(\pi)}^2 \leq C_0 \, n^{-d/2}$$

*for all $n \geq 1$.*

As the proof of Proposition 2.2 is a straightforward adaption of that of Theorem 4.1 in [4], it is also deferred to Appendix A of this paper. A similar result for non-reversible Markov chains on a finite state space can also be found in [5].

We now explain how to establish the upper bound (2.3). Fix $\omega \in \Omega$, and let $Q_\omega$ and $Q_\omega^*$ denote the transition kernel of the random walk under $P_{0,\omega}$ and $P_{0,\omega}^*$, respectively. By (1.4) and our standing Assumptions (b) and (c), the random walk under $P_{0,\omega}$ clearly satisfies the assumptions of Proposition 2.1, the constants $c$, $C$, $\varepsilon$, $N$ and $B$ not depending on $\omega$. Hence, by Propositions 2.1 and 2.2, there exist $m \geq 1$ and $C_0 > 0$ not depending on $\omega$ such that

$$\|(Q_\omega)^{mn}\|_{L^1(\pi_\omega) \to L^2(\pi_\omega)}^2 \leq C_0 \, n^{-d/2}$$

for all $n \geq 1$. Since $Q_\omega$ is a contraction on $L^2(\pi_\omega)$, this implies that

$$\|(Q_\omega)^n\|_{L^1(\pi_\omega) \to L^2(\pi_\omega)}^2 \leq C_0' \, n^{-d/2}$$

for all $n \geq m$, where $C_0' := C_0 \, (m+1)^{d/2}$. Moreover, since the random walk under $P_{0,\omega}^*$ is of the same type as the random walk under $P_{0,\omega}$, we also have

$$\|(Q_\omega^*)^n\|_{L^1(\pi_\omega) \to L^2(\pi_\omega)}^2 \leq C_0' \, n^{-d/2}$$

for all $n \geq m$, possibly after replacing $C_0'$ and $m$ with some larger constants. Therefore, for $n = 2k$ even ($k \geq m$), it follows that

$$\sup_{x,y} \frac{Q_\omega^{2k}(x, y)}{\pi_\omega(\{y\})} = \sup_{x,y} \frac{\langle Q_\omega^{2k} \delta_y, \delta_x \rangle_{\pi_\omega}}{\pi_\omega(\{x\}) \pi_\omega(\{y\})}$$

$$= \sup_{x,y} \frac{\langle (Q_\omega)^k \delta_y, (Q_\omega^*)^k \delta_x \rangle_{\pi_\omega}}{\pi_\omega(\{x\})\pi_\omega(\{y\})}$$

$$\leq \sup_{x,y} \frac{\|(Q_\omega)^k \delta_y\|_{L^2(\pi_\omega)} \|(Q_\omega^*)^k \delta_x\|_{L^2(\pi_\omega)}}{\|\delta_x\|_{L^1(\pi_\omega)} \|\delta_y\|_{L^1(\pi_\omega)}} \leq C_0' k^{-d/2},$$

where the first inequality follows from the Cauchy–Schwarz inequality. Also, for $n = 2k+1$ odd ($k \geq m$), it follows that

$$\sup_{x,y} \frac{Q_\omega^{2k+1}(x,y)}{\pi_\omega(\{y\})} = \sup_{x,y} \sum_z Q_\omega(x,z) \frac{Q_\omega^{2k}(z,y)}{\pi_\omega(\{y\})} \leq C_0' k^{-d/2}.$$

Thus, there is a constant $C_1$ such that (2.3) holds for all $n \geq 2m$, for any $\omega \in \Omega$. By our Assumption (b), replacing $C_1$ with a larger constant if necessary, (2.3) also holds for $n = 1, \ldots, 2m-1$, for any $\omega \in \Omega$. This establishes the desired bound.

## 3. The construction of the corrector function

For each $\omega \in \Omega$, the *chain of environments viewed from the particle* is the $(\Omega, \mathcal{F})$-valued process $(\overline{\omega}_n)_{n \in \mathbb{N}}$ defined by

$$\overline{\omega}_n := T_{X_n}\omega, \quad n \in \mathbb{N}.$$

It is easy to check that for each $\omega \in \Omega$, under the probability measure $P_{0,\omega}$, $(\overline{\omega}_n)_{n \in \mathbb{N}}$ is a Markov chain with transition kernel

$$Rf(\overline{\omega}) := \sum_{z \in \Lambda} f(T_z \overline{\omega}) p_z(\overline{\omega})$$

and initial distribution

$$P_{0,\omega}(\overline{\omega}_0 = \omega) = 1$$

(see e.g. Proposition 1.1 in [3]). Similarly, under the product measure $\mathbb{P} \times P_{0,\cdot}$, $(\overline{\omega}_n)_{n \in \mathbb{N}}$ is a Markov chain with transition kernel $R$ and initial distribution $\mathbb{P}$ (see e.g. Proposition 1.1 in [3]). It is well known that in many interesting "applications," there exists a probability measure $\mathbb{Q}$ on $(\Omega, \mathcal{F})$ which is equivalent to $\mathbb{P}$ and which is invariant for $R$, i.e. we have

$$\int Rf(\omega)\mathbb{Q}(d\omega) = \int f(\omega)\mathbb{Q}(d\omega)$$

for all bounded measurable functions $f$. The existence of such a probability measure $\mathbb{Q}$ is often a prerequisite for the closer investigation of the RWRE.

For a RWRE admitting a bounded cycle representation, such an invariant probability measure $\mathbb{Q}$ for the chain of environment is given by $\mathbb{Q}(d\omega) := Z^{-1} M(\omega) \mathbb{P}(d\omega)$, where $M$ is the positive random variable from the Introduction and $Z := \int_\Omega M(\omega) \mathbb{P}(d\omega)$ is a normalization factor. In the random conductance model considered by Sidoravicius and Sznitman [14] (see Example (a) in the Introduction), the chain of environment is even reversible with respect to $\mathbb{Q}$. However, for a general RWRE admitting a bounded cycle representation, reversibility is usually lost when one of the underlying cycles $C_i$ has length $n_i \geq 3$. At least, we have an explicit description of the time-reversed process: It is induced by the reversed cycles.

Indeed, from the discussion in Section 2, it is clear that the chain of environments associated with the time-reversed RWRE has the transition kernel

$$R^* f(\overline{\omega}) := \sum_{z \in \Lambda^*} f(T_z \overline{\omega}) p_z^*(\overline{\omega}),$$

where $p_z^*(\omega)$ is defined analogously to $p_z(\omega)$ and $\Lambda^*$ is defined analogously to $\Lambda$. To prove that $R^*$ is in fact the adjoint of $R$ in $L^2(\mathbb{Q})$, we have to check that

$$\langle R^*f, g\rangle_\mathbb{Q} = \langle f, Rg\rangle_\mathbb{Q} \tag{3.1}$$

for all non-negative measurable functions $f, g$, where $\langle \cdot, \cdot \rangle_\mathbb{Q}$ denotes the scalar product in $L^2(\mathbb{Q})$. (Note that this also proves our claim that $\mathbb{Q}$ is an invariant measure both for $R$ and for $R^*$.) Now,

$$\langle f, Rg\rangle_\mathbb{Q} = Z^{-1} \cdot \sum_{z \in \Lambda} \int f(\omega) g(T_z\omega) p_z(\omega) M(\omega) \mathbb{P}(d\omega),$$

$$\langle R^*f, g\rangle_\mathbb{Q} = Z^{-1} \cdot \sum_{z \in \Lambda^*} \int f(T_z\omega) g(\omega) p_z^*(\omega) M(\omega) \mathbb{P}(d\omega)$$

$$= Z^{-1} \cdot \sum_{z \in \Lambda^*} \int f(\omega) g(T_{-z}\omega) p_z^*(T_{-z}\omega) M(T_{-z}\omega) \mathbb{P}(d\omega),$$

where the last step uses the translation invariance of $\mathbb{P}$. In view of $\Lambda^* = -\Lambda$ and (2.2), this proves (3.1).

Thus, for a RWRE admitting a bounded cycle representation, there always exists an invariant probability measure $\mathbb{Q} \sim \mathbb{P}$ for the transition kernel $R$. By a straightforward adaption of the proof of Theorem 1.2 in [3] (making use of our standing assumptions (a′) and (c)), this implies that the Markov chain with transition kernel $R$ and initial distribution $\mathbb{Q}$ is ergodic, and there exists at most one invariant probability measure $\mathbb{Q} \sim \mathbb{P}$ for the transition kernel $R$. We will therefore call $\mathbb{Q}$ the invariant probability measure in the sequel.

A quite general approach to deriving invariance principles for RWRE's, which is also used in Sidoravicius and Sznitman [14] and which goes back to [6], is as follows: One constructs a *corrector function* $\chi : \mathbb{Z}^d \times \Omega \to \mathbb{R}^d$ such that for $\mathbb{P}$-a.e. $\omega \in \Omega$, the process $(M_n^\omega)_{n \in \mathbb{N}}$ defined by

$$M_n^\omega := X_n + \chi(X_n, \omega), \quad n \in \mathbb{N},$$

is a martingale under $P_{0,\omega}$. Then one applies the invariance principle for martingales to $(M_n^\omega)_{n \in \mathbb{N}}$, and the (demanding) rest of the proof consists in showing that the contribution of the corrector function is negligible in the limit.

Since the arguments for the construction of the corrector function used in [14] heavily rely on the reversibility of the chain of environments with respect to its invariant distribution $\mathbb{Q}$, they do not apply in the case of a RWRE admitting a bounded cycle representation. We therefore use some different arguments, taken from the field of (asymmetric) exclusion processes (see [9, 10, 15]).

For any $\lambda > 0$, let $u_\lambda$ denote the solution of the *resolvent equation*

$$(\lambda - L)u_\lambda = d_0. \tag{3.2}$$

Here, $L := R - I$ is the (discrete-time) generator of the chain of environments, and $d_0$ is the *local drift at the origin*, which is given by

$$d_0(\omega) := \sum_{z \in \Lambda} z p_z(\omega), \quad \omega \in \Omega.$$

Note that $u_\lambda$ is a well-defined element of $L^2(\mathbb{Q})$, since $d_0$ is an element of $L^2(\mathbb{Q})$ (being bounded) and the operator $\lambda - L$ is invertible in $L^2(\mathbb{Q})$ for any $\lambda > 0$.

Also note that, due to our assumption (b), we have

$$\left(\frac{c}{C}\right) \int f \, d\mathbb{P} \leq \int f \, d\mathbb{Q} \leq \left(\frac{C}{c}\right) \int f \, d\mathbb{P} \tag{3.3}$$

for any measurable function $f \geq 0$. In particular, we have $L^2(\mathbb{Q}) = L^2(\mathbb{P})$, and convergence in $L^2(\mathbb{Q})$ and convergence in $L^2(\mathbb{P})$ are equivalent. Furthermore, since $\mathbb{P}$ is translation invariant, we have

$$\int f \circ T_x \, d\mathbb{Q} \leq \left(\frac{C}{c}\right) \int f \circ T_x \, d\mathbb{P} = \left(\frac{C}{c}\right) \int f \, d\mathbb{P} \leq \left(\frac{C}{c}\right)^2 \int f \, d\mathbb{Q} \tag{3.4}$$

for any measurable function $f \geq 0$ and any $x \in \mathbb{Z}^d$. Thus, $f \in L^2(\mathbb{Q})$ implies $f \circ T_x \in L^2(\mathbb{Q})$ for any $x \in \mathbb{Z}^d$.

**Lemma 3.1.** *For each $x \in \mathbb{Z}^d$ with $\|x\|_2 = 1$, the limit*

$$\lim_{\lambda \to 0} (u_\lambda \circ T_x - u_\lambda)$$

*exists in $L^2(\mathbb{Q})$.*

**Proof.** Since $L^2(\mathbb{Q})$ is complete, it suffices to show that

$$\lim_{\substack{\lambda_1 \to 0 \\ \lambda_2 \to 0}} \|(u_{\lambda_1} \circ T_x - u_{\lambda_1}) - (u_{\lambda_2} \circ T_x - u_{\lambda_2})\|_{L^2(\mathbb{Q})} = 0.$$

Note that the norm can be rewritten as $\|(u_{\lambda_1} - u_{\lambda_2}) \circ T_x - (u_{\lambda_1} - u_{\lambda_2})\|_{L^2(\mathbb{Q})}$. We therefore derive an upper bound for $\|f \circ T_x - f\|_{L^2(\mathbb{Q})}$, where $f \in L^2(\mathbb{Q})$ is arbitrary. By our standing assumption (c) or (c'), for each $\omega \in \Omega$, there exist $1 \leq n(\omega) \leq N$ and (pairwise different) $x_0(\omega), x_1(\omega), \ldots, x_{n(\omega)}(\omega) \in \mathbb{Z}^d$ such that $x_0(\omega) = 0$, $x_{n(\omega)}(\omega) = x$ and $p_{x_i(\omega) - x_{i-1}(\omega)}(T_{x_{i-1}(\omega)}\omega) \geq \varepsilon_0$ for $i = 1, \ldots, n(\omega)$. Thus we obtain

$$((f \circ T_x - f)(\omega))^2 = \left(\sum_{i=1}^{n(\omega)} (f \circ T_{x_i(\omega)} - f \circ T_{x_{i-1}(\omega)})(\omega)\right)^2$$

$$\leq n(\omega) \sum_{i=1}^{n(\omega)} ((f \circ T_{x_i(\omega)} - f \circ T_{x_{i-1}(\omega)})(\omega))^2$$

$$\leq N\varepsilon_0^{-1} \sum_{i=1}^{n(\omega)} ((f \circ T_{x_i(\omega)} - f \circ T_{x_{i-1}(\omega)})(\omega))^2 \cdot p_{x_i(\omega) - x_{i-1}(\omega)}(T_{x_{i-1}(\omega)}\omega)$$

$$\leq N\varepsilon_0^{-1} \sum_{\|z\|_\infty \leq NB} \sum_{z' \in \Lambda} ((f \circ T_{z+z'} - f \circ T_z)(\omega))^2 \cdot p_{z'}(T_z \omega).$$

Here, $\|z\|_\infty := \max_{i=1,\ldots,d} |z_i|$, and $B := \max\{\|z\|_\infty \colon z \in \Lambda\}$. Taking norms and using (3.4), it follows that

$$\|f \circ T_x - f\|_{L^2(\mathbb{Q})}^2 \leq N\varepsilon_0^{-1} \sum_{\|z\|_\infty \leq NB} \int \sum_{z' \in \Lambda} (f \circ T_{z+z'} - f \circ T_z)^2 \cdot p_{z'}(T_z) \, d\mathbb{Q}$$

$$\leq N\varepsilon_0^{-1} (2NB+1)^d \left(\frac{C}{c}\right)^2 \int \sum_{z' \in \Lambda} (f \circ T_{z'} - f)^2 \cdot p_{z'} \, d\mathbb{Q}.$$

Since $\mathbb{Q}$ is an invariant measure for $R$, we have

$$\int \sum_{z \in \Lambda} (f \circ T_z)^2 \cdot p_z \, d\mathbb{Q} = \int Rf^2 \, d\mathbb{Q} = \int f^2 \, d\mathbb{Q} = \int \sum_{z \in \Lambda} f^2 \cdot p_z \, d\mathbb{Q},$$

so that the last integral can be rewritten as

$$\int \sum_{z \in \Lambda} (f \circ T_z - f)^2 \cdot p_z \, d\mathbb{Q} = 2 \int \sum_{z \in \Lambda} f \cdot (f - f \circ T_z) \cdot p_z \, d\mathbb{Q} = 2\langle f, (-L)f \rangle_\mathbb{Q}.$$

It therefore remains to show that

$$\lim_{\substack{\lambda_1 \to 0 \\ \lambda_2 \to 0}} \langle (u_{\lambda_1} - u_{\lambda_2}), (-L)(u_{\lambda_1} - u_{\lambda_2}) \rangle_\mathbb{Q} = 0.$$

By Lemma 2.5.1 in [9], this is true if

$$d_0 \in H_{-1} \quad \text{and} \quad \sup_{\lambda > 0} \|Lu_\lambda\|_{-1} < \infty. \tag{3.5}$$

We refer to Chapters 1 and 2 in [9] for the definitions of the Hilbert spaces $H_{+1}$ and $H_{-1}$ and their respective norms $\|\cdot\|_{+1}$ and $\|\cdot\|_{-1}$. By Sections 2.4 and 2.6 in [9], if the chain of environments satisfies the *sector condition*, i.e. there exists a constant $C_0 \in (0, \infty)$ such that

$$\langle f, (-L)g \rangle_\mathbb{Q}^2 \leq C_0 \cdot \langle f, (-L)f \rangle_\mathbb{Q} \cdot \langle g, (-L)g \rangle_\mathbb{Q}$$

for all $f, g \in L^2(\mathbb{Q})$, then the second condition in (3.5) already follows from the first condition in (3.5). The proof is therefore completed by the subsequent two lemmas. □

**Lemma 3.2.** *Suppose that the RWRE admits a bounded cycle representation. Then the chain of environments satisfies the sector condition, i.e. there exists a constant $C_0 \in (0, \infty)$ such that*

$$\langle f, (-L)g \rangle_\mathbb{Q}^2 \leq C_0 \cdot \langle f, (-L)f \rangle_\mathbb{Q} \cdot \langle g, (-L)g \rangle_\mathbb{Q}$$

*for all $f, g \in L^2(\mathbb{Q})$.*

**Lemma 3.3.** *Suppose that the RWRE admits a bounded cycle representation. Then $d_0 \in H_{-1}$.*

Since the calculations needed for the proof of these lemmas are very similar to those in Section 7.5 in [9] (who treat the special case that the random variable $M$ is constant), they are deferred to Appendix B.

We have just seen that the limits $\lim_{\lambda \to 0}(u_\lambda \circ T_x - u_\lambda)$ ($x \in \mathbb{Z}^d$, $\|x\|_2 = 1$) exist in $L^2(\mathbb{Q})$, and therefore also in $L^2(\mathbb{P})$. Furthermore, since $d_0 \in H_{-1}$ by Lemma 3.3, it follows from resolvent equation (3.2) that $\lim_{\lambda \to 0}(\lambda u_\lambda) = 0$ in $L^2(\mathbb{Q})$ (see Eq. (2.4.3) in [9]), and therefore also in $L^2(\mathbb{P})$. Hence, by considering a suitable subfamily $(\lambda')$ instead of $(\lambda)$, we may assume that the limits

$$\lim_{\lambda' \to 0}(\lambda' u_{\lambda'}(\omega)) = 0$$

and

$$\lim_{\lambda' \to 0}(u_{\lambda'}(T_x\omega) - u_{\lambda'}(\omega)) =: G_x(\omega), \quad x \in \mathbb{Z}^d, \|x\|_2 = 1,$$

exist for $\mathbb{P}$-almost all $\omega \in \Omega$. The random variables $G_x$ thus defined have the following important properties (see also p. 224 in [14]):

**Lemma 3.4.** *The random variables $G_x$ have the following properties:*

(a) *For each $x \in \mathbb{Z}^d$ with $\|x\|_2 = 1$,*

$$\int G_x \, d\mathbb{P} = 0.$$

(b) *If $(x_0, x_1, \ldots, x_n)$ is a sequence in $\mathbb{Z}^d$ such that $\|x_i - x_{i-1}\|_2 = 1$ for all $i = 1, \ldots, n$,*

$$\sum_{i=1}^n G_{x_i - x_{i-1}} \circ T_{x_{i-1}} = \lim_{\lambda' \to 0}(u_{\lambda'} \circ T_{x_n} - u_{\lambda'} \circ T_{x_0}) \quad \mathbb{P}\text{-}a.s.$$

**Proof.** Part (a) follows from the facts that, by translation invariance of $\mathbb{P}$, we have

$$\int (u_{\lambda'} \circ T_x - u_{\lambda'}) \, d\mathbb{P} = 0$$

for all $\lambda' > 0$ and that we have convergence in $L^2(\mathbb{P})$.

Part (b) follows from the facts that we have

$$\sum_{i=1}^{n} (u_{\lambda'} \circ T_{x_i - x_{i-1}} - u_{\lambda'}) \circ T_{x_{i-1}} = u_{\lambda'} \circ T_{x_n} - u_{\lambda'} \circ T_{x_0}$$

for all $\lambda' > 0$ and that we have almost sure convergence. □

We now turn to the construction of the corrector function. For each $x \in \mathbb{Z}^d$ and each $\omega \in \Omega$, set

$$\chi(x,\omega) := \sum_{i=1}^{n} G_{x_i - x_{i-1}} \circ T_{x_{i-1}},$$

where $(x_0, \ldots, x_n)$ is an arbitrary sequence such that $x_0 = 0$, $x_n = x$, and $\|x_i - x_{i-1}\|_2 = 1$ for all $i = 1, \ldots, n$. It follows from Lemma 3.4(b) that $\chi(x, \cdot)$ is a well-defined random variable. The corrector function has the following important properties (see also p. 224 in [14]):

**Lemma 3.5.** *The corrector function has the following properties:*

(a) *For $\mathbb{P}$-almost all $\omega \in \Omega$, $\chi(x+y,\omega) = \chi(x,\omega) + \chi(y, T_x\omega)$ for all $x, y \in \mathbb{Z}^d$.*
(b) *For $\mathbb{P}$-almost all $\omega \in \Omega$, $\sum_{z \in \Lambda} \chi(z,\omega) \, p_z(\omega) = -d_0(\omega)$.*
(c) *For $\mathbb{P}$-almost all $\omega \in \Omega$, the process $M_n^\omega := X_n + \chi(X_n, \omega)$ is a martingale under $P_{0,\omega}$.*

**Proof.** In view of the definition of the corrector function and Lemma 3.4(b), parts (a) and (b) follow from the relations

$$(u_\lambda \circ T_{x+y} - u_\lambda \circ T_0) = (u_\lambda \circ T_x - u_\lambda \circ T_0) + (u_\lambda \circ T_y - u_\lambda \circ T_0) \circ T_x$$

and

$$\sum_{z \in \Lambda} (u_\lambda \circ T_z - u_\lambda \circ T_0) p_z = L u_\lambda = \lambda u_\lambda - d_0.$$

For part (c), first note that, by (a), for $\mathbb{P}$-almost all $\omega \in \Omega$,

$$M_{n+1}^\omega - M_n^\omega = X_{n+1} - X_n + \chi(X_{n+1}, \omega) - \chi(X_n, \omega)$$
$$= X_{n+1} - X_n + \chi(X_{n+1} - X_n, T_{X_n}\omega).$$

Thus, using (b), it follows that for $\mathbb{P}$-almost all $\omega \in \Omega$,

$$E_{0,\omega}(M_{n+1}^\omega - M_n^\omega | X_0, \ldots, X_n) = \sum_{z \in \Lambda} (z + \chi(z, T_{X_n}\omega)) p_z(T_{X_n}\omega)$$
$$= \sum_{z \in \Lambda} z p_z(T_{X_n}\omega) + \sum_{z \in \Lambda} \chi(z, T_{X_n}\omega) p_z(T_{X_n}\omega)$$
$$= d_0(T_{X_n}\omega) - d_0(T_{X_n}\omega)$$
$$= 0.$$

□

## 4. Proof of the main theorem

A detailed analysis of the proof of Theorem 1.1 in [14] reveals that it does not use any special properties of the random conductance model (in particular, it does not use reversibility), but only

- the properties (1.10)–(1.14) concerning the corrector function,
- the upper Gaussian bound (1.16) for the transition kernel,
- the Markov property and the bounded range of the random walk,
- the Markov property and the ergodicity of the chain of environments,
- the assumption (a).

(See also Remark 1.3 in [14].) Therefore, since we have seen in the preceding sections that these statements remain true for a RWRE admitting a bounded cycle representation, it follows that the proof also applies to this model.

In particular, the limit distribution is also a $d$-dimensional Brownian motion with mean zero and covariance matrix

$$A = \left( \int \sum_{z \in \Lambda} p_z(\omega) \langle z + \chi(z,\omega), e_i \rangle \langle z + \chi(z,\omega), e_j \rangle \mathbb{Q}(d\omega) \right)_{i,j} \quad (4.1)$$

(see Eq. (1.15) in [14]). Here $\langle \cdot, \cdot \rangle$ denotes the scalar product in $\mathbb{R}^d$. Naturally, in contrast to the random conductance model with i.i.d. couplings, the covariance matrix need not be of the form $\sigma^2 I$ with $\sigma^2 \geq 0$ anymore. The non-degeneracy of the covariance matrix follows from a similar argument as in Remark 1.2 in [14].

Indeed, let $x \in \mathbb{R}^d \setminus \{0\}$. Then there exists a vector $e \in \mathbb{Z}^d$ with $\|e\|_2 = 1$ such that $\langle e, x \rangle > 0$. By Lemma 3.4(a), it follows that

$$\int \langle e + \chi(e,\omega), x \rangle \mathbb{P}(d\omega) > 0,$$

so that

$$\mathbb{P}(\{\omega: \langle e + \chi(e,\omega), x \rangle > 0\}) > 0.$$

By our standing assumption (c) or (c'), it further follows that there exist $n \in \mathbb{N}$ and $z_0, \ldots, z_n \in \mathbb{Z}^d$ such that $z_0 = 0$, $z_n = e$ and

$$\mathbb{P}(\{\omega: \langle e + \chi(e,\omega), x \rangle > 0 \text{ and } p_{z_i - z_{i-1}}(T_{z_{i-1}}\omega) > \varepsilon_0 \ \forall i = 1, \ldots, n\}) > 0.$$

Since

$$e + \chi(e,\omega) = \sum_{i=1}^{n} (z_i - z_{i-1}) + \chi(z_i - z_{i-1}, T_{z_{i-1}}\omega)$$

by Lemma 3.5(a), it follows that

$$\mathbb{P}(\{\omega: \langle z_i - z_{i-1} + \chi(z_i - z_{i-1}, T_{z_{i-1}}\omega), x \rangle > 0 \text{ and } p_{z_i - z_{i-1}}(T_{z_{i-1}}\omega) > \varepsilon_0\}) > 0$$

for some $i = 1, \ldots, n$. By translation invariance of $\mathbb{P}$, this reduces to

$$\mathbb{P}(\{\omega: \langle z + \chi(z,\omega), x \rangle > 0 \text{ and } p_z(\omega) > \varepsilon_0\}) > 0$$

for some $z \in \Lambda$. We may therefore conclude that

$$x^T A x = \int \sum_{z \in \Lambda} p_z(\omega) |\langle z + \chi(z,\omega), x \rangle|^2 \mathbb{Q}(d\omega) > 0,$$

which proves the non-degeneracy of $A$.

## Appendix A. The Nash inequality

This section is devoted to the proofs of Propositions 2.1 and 2.2 from Section 2. We will find it convenient to work with the maximum norm $\|\cdot\|_\infty$ instead of the Euclidean norm $\|\cdot\|_2$.

Our results are based on the following assumption on a transition kernel $Q$ on $\mathbb{Z}^d$:

**Assumption A.1.** *There exist $K \in \mathbb{N}$ and $\delta > 0$ such that for all $x, y \in \mathbb{Z}^d$ with $\|y - x\|_\infty \leq 3K + 1$, there exists $y' \in \mathbb{Z}^d$ with $\|y' - y\|_\infty \leq K$ and $Q(x, y') \geq \delta$.*

Proposition 2.1 will follow immediately from the next two lemmas:

**Lemma A.2.** *Let $(X_n)_{n \in \mathbb{N}}$ be a random walk on $\mathbb{Z}^d$ with transition kernel $Q$ and invariant measure $\pi$ for which there exist constants $0 < c \leq C < \infty$ with $c \leq \pi(\{x\}) \leq C$ for all $x \in \mathbb{Z}^d$. Suppose that the transition kernel $Q$ is strongly irreducible, i.e.*

$$\exists \varepsilon > 0, \ \exists N \in \mathbb{N}, \ \forall x \in \mathbb{Z}^d, \ \forall e \in \mathbb{Z}^d \colon \|e\|_2 = 1, \ \exists n \leq N, \quad Q^n(x, x + e) > \varepsilon,$$

*and of bounded range, i.e.*

$$\exists B > 0, \ \forall x, y \in \mathbb{Z}^d, \quad \|y - x\|_\infty > B \implies Q(x, y) = 0.$$

*Then there exists an $m \geq 1$ (depending only on $c$, $C$, $\varepsilon$, $N$ and $B$) such that the transition kernel $(Q^m)^*(Q^m)$ satisfies Assumption A.1.*

**Proof.** First of all, observe that the adjoint $Q^*$ of $Q$ with respect to $\pi$ is also strongly irreducible and of finite range, since

$$(Q^*)^n(x, y) = \frac{\pi(\{y\})}{\pi(\{x\})} Q^n(y, x) \in \left[ \frac{c}{C} Q^n(y, x); \frac{C}{c} Q^n(y, x) \right]$$

for all $n \geq 1$, $x \in \mathbb{Z}^d$, $y \in \mathbb{Z}^d$. Hence, replacing $\varepsilon$ by some smaller constant $\varepsilon'$ if necessary, we may assume that the assumptions of the lemma are satisfied both for $Q$ and for $Q^*$.

It follows from the assumptions that for all $x \in \mathbb{Z}^d$ and for all $y \in \mathbb{Z}^d$ with $\|y\|_2 = 1$, there exist $n \leq N$ and a sequence $(x_0, \ldots, x_n)$ in $\mathbb{Z}^d$ such that $x_0 = x$, $x_n = x + y$ and $Q(x_{i-1}, x_i) \geq \varepsilon_0 := \varepsilon/(2B+1)^{dN}$ for all $i = 1, \ldots, n$. Such a sequence will also be called an $\varepsilon_0$-path from $x$ to $x + y$. By replacing $N$ with $2N$, we may also assume that for each $x \in \mathbb{Z}^d$, there exists an $\varepsilon_0$-cycle for $x$ (i.e. an $\varepsilon_0$-path from $x$ to $x$).

Pick $K := NB$ and $L := 3K + 1$. Then, for all $x, y \in \mathbb{Z}^d$ with $\|y - x\|_\infty \leq L$, there is an $\varepsilon_0$-path from $x$ to $y$ of length $\leq dLN$. (Indeed, there is certainly a nearest-neighbor path from $x$ to $y$ of length $\leq dL$, and by strong irreducibility, each step of this path can be replaced with at most $N$ steps of the random walk associated with $Q$.) As $\varepsilon_0$-paths can be extended by adding $\varepsilon_0$-cycles, there is also an $\varepsilon_0$-path from $x$ to $x+y$ of length $\in \{dLN + 1, \ldots, dLN + N\}$. If we cut the path immediately after step $dLN$ (at site $y'$, say) we have $\|y' - y\|_\infty \leq NB$ (because we get to site $y$ exactly after at most $N$ additional steps, each of which has size $\leq B$) and $Q^{dLN}(x, y') \geq \varepsilon_0^{dLN}$. The same argument applied to $Q^*$ yields a site $y''$ such that $\|y'' - y\|_\infty \leq NB$ and $(Q^*)^{dLN}(y', y'') \geq \varepsilon_0^{dLN}$. Summarizing, we have $\|y'' - y\|_\infty \leq NB$ and $Q^{dLN}(Q^{dLN})^*(x, y'') \geq \varepsilon_0^{2dLN}$. Putting $m := dLN$ and $\delta := \varepsilon_0^{2dLN}$ and exchanging the roles of $Q$ and $Q^*$ completes the proof. $\square$

To see the connection to the next lemma, note that the transition kernel $(Q^m)^*Q^m$ constructed above has $\pi$ as an invariant and reversible measure.

**Lemma A.3.** *Suppose that $Q$ is a transition kernel satisfying Assumption A.1 and that $\pi$ is an invariant and reversible measure for $Q$ such that there exist constants $0 < c \leq C < \infty$ with $c \leq \pi(\{x\}) \leq C$ for all $x \in \mathbb{Z}^d$. Then the transition kernel $Q$ satisfies the $d$-dimensional isoperimetric inequality*

$$\exists \kappa > 0, \ \forall A \subset_{\text{finite}} \mathbb{Z}^d, \quad (\pi(A))^{1-1/d} \leq \kappa a(\partial A)$$

*and therefore the d-dimensional Nash inequality*

$$\exists \kappa' > 0, \ \forall f \in \ell_0(\mathbb{Z}^d), \quad \mathcal{E}_Q(f,f) \geq \kappa' \|f\|_{L^2(\pi)}^{2+4/d} \|f\|_{L^1(\pi)}^{-4/d}.$$

*More precisely, $\kappa$ and $\kappa'$ are constants depending only on $c$ and $C$ as well as on $K$ and $\delta$ (from Assumption A.1).*

Here, we set $\pi(A) := \sum_{x \in A} \pi(\{x\})$, $a(x,y) := \pi(x)Q(x,y) = \pi(y)Q(y,x)$, $\partial A :=$ the set of edges having one endpoint in $A$ and one endpoint in $\mathbb{Z}^d \setminus A$, $a(\partial A) := \sum_{e \in \partial A} a(e)$, $\ell_0(\mathbb{Z}^d) :=$ the set of functions on $\mathbb{Z}^d$ with finite support. See Sections 4 and 14 in [16] for further details.

**Proof of Lemma A.3.** By Proposition 14.1 in [16], it suffices to establish the $d$-dimensional isoperimetric inequality. To this end, we compare the random walk associated with $Q$ (or, more precisely, some kind of renormalization of it) to the standard random walk on $\mathbb{Z}^d$, for which the $d$-dimensional isoperimetric inequality is known to hold (see p. 45 in [16]).

Take $L := 2K + 1$ (where $K$ is the constant from Assumption A.1), and partition $\mathbb{Z}^d$ into $d$-dimensional cubes of length $L$. The cube containing $x$ will be denoted by $C(x)$, and the set of all cubes will be denoted by $C(\mathbb{Z}^d)$. Write $C(x) \sim C(y)$ if $C(x)$ and $C(y)$ have a common "face," and observe that the graph $C(\mathbb{Z}^d)$ thus obtained is isomorphic to $\mathbb{Z}^d$ with the nearest-neighbor topology.

Let $\overline{\pi}$ denote the uniform measure on $C(\mathbb{Z}^d)$, let $\overline{Q}$ be the transition kernel on $C(\mathbb{Z}^d)$ corresponding to the standard random walk on $\mathbb{Z}^d$, and let $\overline{a}$ denote the associated conductance. Then the $d$-dimensional isoperimetric inequality for the standard random walk on $\mathbb{Z}^d$ (see above) states that

$$\exists \kappa > 0, \ \forall \overline{A} \subset_{\text{finite}} C(\mathbb{Z}^d), \quad (\overline{\pi}(\overline{A}))^{1-1/d} \leq \kappa \overline{a}(\partial \overline{A}).$$

Now consider $A \subset_{\text{finite}} \mathbb{Z}^d$ for the random walk with transition kernel $Q$. Then $\overline{A} := \{C(x) \in C(\mathbb{Z}^d) : x \in A\}$ is a finite subset of $C(\mathbb{Z}^d)$, and clearly

$$\pi(A) \leq C \cdot L^d \cdot \overline{\pi}(\overline{A}),$$

$C$ denoting the upper bound on $\pi(\{x\})$, $x \in \mathbb{Z}^d$.

Now let $\overline{e} = (\overline{x}, \overline{y}) \in \partial \overline{A}$, where $\overline{x} \in \overline{A}$ and $\overline{y} \in C(\mathbb{Z}^d) \setminus \overline{A}$, say. Let $x \in \overline{x} \cap A$, and pick $y \in \overline{y} \subset \mathbb{Z}^d \setminus A$ such that $Q(x,y) \geq \delta$. Such a vertex $y$ exists by our assumption on $Q$. Indeed, since $x$ has $\|\cdot\|_\infty$-distance $\leq K$ from the center of $\overline{x}$ and therefore $\|\cdot\|_\infty$-distance $\leq 3K + 1$ from the center of $\overline{y}$, our assumption on $Q$ ensures the existence of an element $y \in \mathbb{Z}^d$ with $\|\cdot\|_\infty$-distance $\leq K$ from the center of $\overline{y}$ (i.e. $y \in \overline{y}$) and $Q(x,y) \geq \delta$. Applying this argument to any $\overline{e} \in \partial \overline{A}$, it follows that

$$c \cdot \delta \cdot \overline{a}(\partial \overline{A}) \leq a(\partial A),$$

$c$ denoting the lower bound on $\pi(\{x\})$, $x \in \mathbb{Z}^d$.

Putting the preceding inequalities together, it follows that the random walk associated with $Q$ satisfies the $d$-dimensional isoperimetric inequality. □

**Remark A.4.** *It follows from the preceding results that for a random walk as in Lemma A.2, the $d$-dimensional Nash inequality holds for the transition kernel $(Q^m)^* Q^m$, for a suitable $m \geq 1$. We mention without proof that one can easily construct (non-reversible and non-translation-invariant) random walks for which the assumptions of Lemma A.2 are satisfied, but for which the $d$-dimensional Nash inequality does not hold for the transition kernel $Q^* Q$.*

We now turn to the proof of Proposition 2.2:

**Proof of Proposition 2.2.** The proof is a straightforward adaption of the proof of Theorem 4.1 in [4].

It is easy to see that if the Nash inequality holds for all $f \in \ell_0(\mathbb{Z}^d)$, then it also holds for all $f \in L^2(\pi)$. Thus, we have

$$\begin{aligned}\|Q^{n+1}f\|^2_{L^2(\pi)} - \|Q^n f\|^2_{L^2(\pi)} &= \langle Q^n f, Q^*QQ^n f\rangle_\pi - \langle Q^n f, Q^n f\rangle_\pi \\ &= -\mathcal{E}_{Q^*Q}(Q^n f, Q^n f) \\ &\leq -\kappa \|Q^n f\|^{-4/d}_{L^1(\pi)} \|Q^n f\|^{2+4/d}_{L^2(\pi)}\end{aligned}$$

for all $n \in \mathbb{N}$. Hence, if $f \in L^1_+(\pi)$ with $\|f\|_{L_1(\pi)} = 1$ and $u_n := \|Q^n f\|^2_{L_2(\pi)}$, $n \in \mathbb{N}$, we are led to the difference equation/inequality

$$u_{n+1} \leq u_n - \kappa u_n^{1+2/d} = u_n(1 - \kappa u_n^{2/d}), \quad n \in \mathbb{N}.$$

We must show that this implies

$$u_n \leq C_0/n^{d/2}, \quad n = 1, 2, \ldots, \tag{A.1}$$

for some $C_0 \in (0, \infty)$.

We can clearly find $C_0 > 0$ (depending only on $c$, $C$ and $\kappa$) such that

$$u_1 \leq C_0 \quad \text{and} \quad \left(1 - \frac{\kappa C_0^{2/d}}{n+1}\right) \leq \left(\frac{n}{n+1}\right)^{d/2} \quad \forall n \geq 1.$$

Then it follows by induction that (A.1) holds. Indeed, suppose that $u_n \leq C_0/n^{d/2}$ for some $n \geq 1$. If $u_n \leq C_0/(n+1)^{d/2}$, then $u_{n+1} \leq C_0/(n+1)^{d/2}$. If $u_n > C_0/(n+1)^{d/2}$, then

$$u_{n+1} \leq u_n(1 - \kappa u_n^{2/d}) < u_n\left(1 - \frac{\kappa C_0^{2/d}}{n+1}\right) \leq \left(\frac{C_0}{n^{d/2}}\right) \cdot \left(\frac{n}{n+1}\right)^{d/2} = \frac{C_0}{(n+1)^{d/2}}.$$

This completes the proof. $\square$

## Appendix B. Some lengthy calculations

In this section we prove Lemmas 3.2 and 3.3. The calculations are almost the same as in Section 7.5 in [9], who treat the special case that the random variable $M$ is constant.

**Proof of Lemma 3.2.** It easily follows from our definitions that for all $\omega \in \Omega$ and all $z \in \Lambda$,

$$p_z(\omega) = \frac{1}{M(\omega)} \sum_{i=1}^{K} \sum_{j=1}^{n_i} W_i(T_{-z_{i,j-1}}\omega) \cdot \mathbf{1}_{\{z_{i,j} - z_{i,j-1} = z\}},$$

where $C_i = (z_{i,0}, \ldots, z_{i,n_i})$. Hence, the generator for the chain of environments has the representation

$$\begin{aligned}Lf(\omega) &= \sum_{z \in \Lambda}(f(T_z \omega) - f(\omega)) \cdot p_z(\omega) \\ &= \frac{1}{M(\omega)} \sum_{i=1}^{K} \sum_{j=1}^{n_i} (f(T_{z_{i,j}-z_{i,j-1}}\omega) - f(\omega)) \cdot W_i(T_{-z_{i,j-1}}\omega),\end{aligned}$$

and

$$\langle f, (-L)g\rangle_\mathbb{Q}$$

$$= -\int \frac{1}{M(\omega)} \sum_{i=1}^{K} \sum_{j=1}^{n_i} f(\omega)(g(T_{z_{i,j}-z_{i,j-1}}\omega) - g(\omega)) \cdot W_i(T_{-z_{i,j-1}}\omega) \, d\mathbb{Q}(\omega)$$

$$= -\int \frac{1}{M(\omega)} \sum_{i=1}^{K} \sum_{j=1}^{n_i} f(T_{z_{i,j-1}}\omega)(g(T_{z_{i,j}}\omega) - g(T_{z_{i,j-1}}\omega)) \cdot W_i(\omega) \, d\mathbb{Q}(\omega).$$

Here we have used the identity

$$\int \frac{f(\omega)}{M(\omega)} \, d\mathbb{Q}(\omega) = \int \frac{f(T_z\omega)}{M(\omega)} \, d\mathbb{Q}(\omega),$$

which follows the translation invariance of $\mathbb{P}$.

In the special case $f = g$, we have

$$-\sum_{j=1}^{n_i} f(T_{z_{i,j-1}}\omega)(f(T_{z_{i,j}}\omega) - f(T_{z_{i,j-1}}\omega)) = \frac{1}{2} \sum_{j=1}^{n_i} (f(T_{z_{i,j}}\omega) - f(T_{z_{i,j-1}}\omega))^2$$

for all $i = 1, \ldots, K$, since $\{z_{i,1}, \ldots, z_{i,n_i}\} = \{z_{i,0}, \ldots, z_{i,n_i-1}\}$, and therefore

$$\langle f, (-L)f \rangle_{\mathbb{Q}} = \frac{1}{2} \int \frac{1}{M(\omega)} \sum_{i=1}^{K} \sum_{j=1}^{n_i} (f(T_{z_{i,j}}\omega) - f(T_{z_{i,j-1}}\omega))^2 \cdot W_i(\omega) \, d\mathbb{Q}(\omega).$$

In the general case, we have

$$\sum_{j=1}^{n_i} f(T_{z_{i,j-1}}\omega)(g(T_{z_{i,j}}\omega) - g(T_{z_{i,j-1}}\omega))$$

$$= \sum_{j=1}^{n_i} (f(T_{z_{i,j-1}}\omega) - f(T_{z_{i,0}}\omega))(g(T_{z_{i,j}}\omega) - g(T_{z_{i,j-1}}\omega))$$

for all $i = 1, \ldots, K$, since $\{z_{i,1}, \ldots, z_{i,n_i}\} = \{z_{i,0}, \ldots, z_{i,n_i-1}\}$, and therefore, by the Cauchy–Schwarz inequality,

$$|\langle f, (-L)g \rangle_{\mathbb{Q}}|^2 \leq \int \frac{1}{M(\omega)} \sum_{i=1}^{K} \sum_{j=1}^{n_i} (f(T_{z_{i,j-1}}\omega) - f(T_{z_{i,0}}\omega))^2 \cdot W_i(\omega) \, d\mathbb{Q}(\omega)$$

$$\times \int \frac{1}{M(\omega)} \sum_{i=1}^{K} \sum_{j=1}^{n_i} (g(T_{z_{i,j}}\omega) - g(T_{z_{i,j-1}}\omega))^2 \cdot W_i(\omega) \, d\mathbb{Q}(\omega).$$

The second factor equals $2 \cdot \langle g, (-L)g \rangle_{\mathbb{Q}}$. Moreover, since

$$\sum_{j=1}^{n_i} (f(T_{z_{i,j-1}}\omega) - f(T_{z_{i,0}}\omega))^2 = \sum_{j=1}^{n_i} \left( \sum_{k=1}^{j-1} (f(T_{z_{i,k}}\omega) - f(T_{z_{i,k-1}}\omega)) \right)^2$$

$$\leq \sum_{j=1}^{n_i} (j-1) \cdot \sum_{k=1}^{j-1} (f(T_{z_{i,k}}\omega) - f(T_{z_{i,k-1}}\omega))^2$$

$$\leq n_i^2 \cdot \sum_{j=1}^{n_i} (f(T_{z_{i,j}}\omega) - f(T_{z_{i,j-1}}\omega))^2$$

for all $i=1,\ldots,K$, the first factor is bounded above by $2 \cdot \max_{i=1,\ldots,K} n_i^2 \cdot \langle f, (-L)f \rangle_{\mathbb{Q}}$. It follows that

$$|\langle f, (-L)g \rangle_{\mathbb{Q}}|^2 \leq 4 \cdot \max_{i=1,\ldots,n} n_i^2 \cdot \langle f, (-L)f \rangle_{\mathbb{Q}} \cdot \langle g, (-L)g \rangle_{\mathbb{Q}}.$$

This proves the lemma. □

**Proof of Lemma 3.3.** It follows from the definitions that a function $V \in L^2(\mathbb{Q})$ belongs to $H_{-1}$ if and only if there exists a constant $C_0 \in (0, \infty)$ such that $|\langle V, f \rangle_{\mathbb{Q}}|^2 \leq C_0 \langle f, (-L)f \rangle_{\mathbb{Q}}$ for all $f \in L^2(\mathbb{Q}) \cap H_{+1}$. (See Section 2.7 in [9].) But now, by similar arguments as in the preceding proof, we have

$$\sum_{z \in \Lambda} z p_z(\omega) = \frac{1}{M(\omega)} \sum_{i=1}^{K} \sum_{j=1}^{n_i} (z_{i,j} - z_{i,j-1}) W_i(T_{-z_{i,j-1}}\omega)$$

$$= \frac{1}{M(\omega)} \sum_{i=1}^{K} \sum_{j=1}^{n_i} z_{i,j}(W_i(T_{-z_{i,j-1}}\omega) - W_i(T_{-z_{i,j}}\omega))$$

and therefore, denoting the scalar product and the Euclidean norm in $\mathbb{R}^d$ by $\langle \cdot, \cdot \rangle$ and $\|\cdot\|_2$,

$$|\langle d_0, f \rangle_{\mathbb{Q}}|^2 = \left( \int \frac{1}{M(\omega)} \sum_{i=1}^{K} \sum_{j=1}^{n_i} \langle z_{i,j}, f(\omega) \rangle \cdot (W_i(T_{-z_{i,j-1}}\omega) - W_i(T_{-z_{i,j}}\omega)) \, \mathrm{d}\mathbb{Q}(\omega) \right)^2$$

$$= \left( \int \frac{1}{M(\omega)} \sum_{i=1}^{K} \sum_{j=1}^{n_i} \langle z_{i,j}, f(T_{z_{i,j}}\omega) - f(T_{z_{i,j-1}}\omega) \rangle \cdot W_i(\omega) \, \mathrm{d}\mathbb{Q}(\omega) \right)^2$$

$$\leq \int \frac{1}{M(\omega)} \sum_{i=1}^{K} \sum_{j=1}^{n_i} \|z_{i,j}\|_2^2 \cdot W_i(\omega) \, \mathrm{d}\mathbb{Q}(\omega)$$

$$\quad \times \int \frac{1}{M(\omega)} \sum_{i=1}^{K} \sum_{j=1}^{n_i} \|f(T_{z_{i,j}}\omega) - f(T_{z_{i,j-1}}\omega)\|_2^2 \cdot W_i(\omega) \, \mathrm{d}\mathbb{Q}(\omega)$$

$$\leq 2 \cdot \max_{i=1,\ldots,K} \sum_{j=1}^{n_i} \|z_{i,j}\|_2^2 \cdot \langle f, (-L)f \rangle_{\mathbb{Q}}.$$

This proves the lemma. □